\documentclass[leqno,12pt]{article}
\usepackage{amsthm,amsfonts,amssymb,amsmath}
\usepackage[all]{xy}
\usepackage{xr-hyper}
\usepackage[colorlinks=true]{hyperref}
\usepackage{graphicx}

\newcommand{\BC}{{\mathbb {C}}}

\newcommand{\BQ}{{\mathbb {Q}}}

\newcommand{\BZ}{{\mathbb {Z}}}

\newcommand{\CL}{{\mathcal {L}}}
\newcommand{\CM}{{\mathcal {M}}}

\newcommand{\CP}{{\mathcal {P}}}
\newcommand{\CQ}{{\mathcal {Q}}}

\newcommand{\RN}{{\mathrm {N}}}

\newcommand{\Alb}{{\mathrm{Alb}}}

\newcommand{\Ch}{{\mathrm{Ch}}}

\newcommand{\Hom}{{\mathrm{Hom}}}

\newcommand{\NS}{{\mathrm{NS}}}

\newcommand{\Pic}{\mathrm{Pic}}

\newcommand{\pair}[1]{\langle {#1} \rangle}

\newcommand{\lra}{\longrightarrow}

\theoremstyle{plain}
\newtheorem{theorem}{Theorem}

\theoremstyle{definition}

\theoremstyle{remark}

\newtheorem{remarks}[theorem]{Remarks}

\numberwithin{equation}{subsection}

\title{Positivity of heights of  codimension 2 cycles\\
over function field of characteristic $0$}
\author{Shou-Wu Zhang\\
Department of Mathematics\\
Columbia University\\
New York, NY 10027}
\begin{document}
\maketitle

In this note, we show how the classical  Hodge index theorem implies the 
Hodge index conjecture of Beilinson \cite{Bei} 
for height pairing of homologically trivial 
codimension $2$ cycles over function field of characteristic $0$. 
Such an index conjecture has been used in our previous paper \cite{Zh}
to deduce the Bogomolov conjecture and a lower bound for Hodge class (or Faltings height)
from some conjectures about metrized graphs which have just been recently proved by 
Zubeyir Cinkir \cite{Cin}. We would like to thank Walter Gubler for pointing out
the lack of  a proof or reference of such an  index theorem in \cite{Zh},
and to Walter Gubler and Klaus K\"unnemann for help to prepare this note.

Let $f: X\lra B$ be a flat morphism between smooth, projective, and
geometrically connected varieties  over a field $k$ of
characteristic $0$. Assume that $B$ has dimension $1$ and that $X$ has
dimension  $n\ge 4$. Let  $\CL$ be an ample line bundle on $X$. The
purpose of this note is to prove the following Hodge index theorem
for heights of codimension $2$ cycles:
\begin{theorem} Let $z\in
\Ch ^2(X)$ be a codimension $2$ cycle such that the following two
conditions hold:
\begin{enumerate}
\item $z\cdot v=0$ for any subvariety
$v$ of dimension $2$ included in a fiber of $f$;
\item $z_K\cdot c_1(\CL_K)^{n-3}=0$ in $\Alb(X_K)$.
\end{enumerate}
Then
$$z^2c_1(\CL)^{n-4}\ge  0.$$
Moreover the equality holds if and only if  $z$ is numerically
equivalent to the  fiber  $u$ of a horizontal divisor $U$ of $X$
over a point of $B$.
\end{theorem}

By the Lefschetz principle, we may assume that $k=\BC$. We will
deduce the theorem from the usual Hodge index theorem on $X$. More
precisely we will find  a horizontal divisor $U$ with rational coefficients
such that for any fiber $u$ over $B$, the class of
 $(z-u)\cdot c_1(\CL)^{n-3}$ vanishes in  $H^{2n-2}(X, \BQ)$.
Then by the Hodge index theorem,
$$z^2=(z-u)^2c_1(\CL)^{n-4}\ge 0.$$
Moreover, if the above intersection is $0$ then $z-u$ is numerically
trivial.

Using cup product with $H^2(X, \BC)$ and Hodge decomposition, the
condition here can be replaced by
$$(z-u)\cdot c_1(\CL)^{n-3}\cdot \alpha =0$$
for any $\alpha \in H^{1,1}(X, \BC)\cap H^2(X, \BZ)$. Since such an
$\alpha$ is a linear combination of the first Chern class of line
bundles, we have an equivalent condition
$$(z-u)\cdot c_1(\CL)^{n-3}\cdot c_1(\CM)=0$$
for any line bundle $\CM$ of $X$.

Let $\ell $ denote the functional on $\Pic (X)$ defined by
$$\ell (\CM):=z\cdot c_1(\CL)^{n-3}\cdot c_1(\CM).$$
Since $z$ is perpendicular to vertical $2$ cycles, this functional
depends only on the image of $\CM$ in $\Pic (X_K)$. Thus we may
consider $\ell$ as a functional on $\Pic (X_K)$.

First we notice that $\ell$ vanishes on $\Pic ^0(X_K)$ since it is
the Neron--Tate height pairing between  the class of $z_K\cdot
c_1(\CL_K)^{n-3}$ in $\Alb (X_K)$ and points on $\Pic ^0(X_K)$.
Since  lacking of reference, we provide a proof of  this fact here. Let
$J_K$ be the Picard variety of $X_K$. Thus $J_K(K)=\Pic ^0(X_K)$.
Let $A_K$ be the dual variety of $J_K$ which is also the Albanese
variety of $X_K$. Let $\CP_K$ be the universal bundle on $A_K\times
J_K$ which is trivial when restrict on $\{0\}\times J_K$ and
$A_K\times \{0\}$. Let $J$ and $A$ be the Neron models of $J$ and
$A$ over $B$ respectively. Then $A\times _B J$ is the Neron model of
$A_K\times J_K$ over $B$. Let $\CP$ be the canonical extension of
$\CP_K$ on $A\times _BJ$, see Moret--Bailly \cite{MB}, Chapter II and III.
 Then for a $B$- point $z\in \Hom (B,
A\times_BJ)$ with components $x$ and $y$ on $A$ and $J$
respectively,
$$\deg z^*\CP=\deg x^*\CP_y=-\pair{x, y}_{NT}.$$
Here $\CP_y$ is the line bundle on $A$ corresponding to $y$ as the
connected component of the Picard variety of $A_K$ is still $J_K$.

After a resolution of singularity $X'\lra X_{B'}$ for a base change
by a finite morphism $B'\lra B$, we may assume that $X$ has a
section $e$ over $B$, and that $z\cdot c_1(\CL)^{n-3}$ is
represented by a linear combination $\sum _i n_i s_i+\sum _j m_j
C_j$ of sections $s_i$ and vertical curves $C_i$ in $X\lra B$. Then
we have an Abel--Jacobi $X_K\lra A_K$ which sends $e_K$ to zero. Let
$\CQ_K$ be the pull-back  of $\CP_K$ on $X_K\times J_K$. Then $\CQ$
is the universal bundle with a trivialization on $\{e_K\}\times
J_K$. Using the universal property of Neron model, the Abel--Jacobi
morphism can be extended to a morphism $X^0\lra A$ where $X^0$ is
the complement of $X$ of the locus $S$ of the singular points on the
fibers. Let $\CQ$ be the pull-back of $\CP$ on $X^0\times _B J$.
Since $S$ has codimension at least two on regular variety $X\times
_B A$, $\CQ$ has a unique extension to $X\times _BA$ which we still
denote as $\CQ$. Now the functional $\ell$ on $\Pic ^0(X_K)=J(B)$
has the following expression: for any $t\in J(B)$,
$$\ell (t)=\sum _i n_i\deg _{s_i}\CQ_t+\sum _j m_j\deg
_{C_j}\CQ_t.$$

The first sum vanishes as it is minus the Neron--Tate height pairing
of $\sum_in_is_{iK}=0\in \Alb (X_K)$ and $t$. For the second  sum, we
let $b_j\in B$ be the image of $C_j$, and $A_j$ be the fiber of $A$
over $b_j$, and $\CQ_j$ be the restriction on $C_j\times A_j$. Then
$$\deg _{C_j}\CQ_t=\deg _{C_j}\CQ_{j, t_{b_j}}.$$
As a function on $t_{b_j}\in A_j$, the right hand side of the above
identity  is locally constant. Thus $\deg _{C_j}\CQ_t$ only takes
finitely many values. Thus, we have shown that $\ell$ takes only
finitely many values on $\Pic ^0(X_K)$. Since $\ell$ is additive, it
must vanishes.

In summary, we have shown that the value  $\ell(\CM)$ depends only
on the image of $\CM$ in the Neron-Sevri group $\NS(X_K)_\BQ=\Pic
(X_K)/\Pic^0(X_K)\otimes \BQ$ with rational coefficient. 
By non-degeneracy of pairing between $\NS(X_K)_\BQ$
and the group $\RN _1(X)_\BQ$ of 1-cycles with rational coefficients on $X_K$
modulo numerical equivalence, we have a dimension 1-cycle $W_K$ with
rational coefficients on $X_K$ such that
$$\ell (\CM)=c_1(\CM_K)\cdot W_K.$$
Using hard Lefschetz on cohomology groups and Hodge conjecture for divisors, 
we have a divisor  $U_K$ with rational
coefficients such that $W_K=U_k\cdot c_1(\CL_K)^{n-3}$ modulo
$\Pic^0(X_K)$. Let $U$ be the Zariski closure of $U_K$ on $X$ and
let $u$ be the fiber over $U$ over some point of $B$ consider as a
subvariety of $X$. Then by the flatness of $U\lra B$, we have
$$\ell (\CM)=c_1(\CM_K)\cdot W_K=c_1(\CM)\cdot c_1(\CL)^{n-3}\cdot u.$$
In other words, we have shown the required property: $$(z-u)\cdot
c_1(\CL)^{n-3}\cdot c_1(\CM)=0, \qquad \forall \CM\in \Pic (X).$$

\begin{remarks}The same proof holds for varieties over number field
or function field with positive characteristic provided the Hodge
index conjecture of Grothendieck and Gillet--Soule for codimension
$2$ cycles.
\end{remarks}

\end{document}